%% file: prime_TCS_R1.tex
\PassOptionsToPackage{dvipsnames}{xcolor} 
\documentclass[11pt,letter]{article}
\usepackage{amssymb,amsmath,amsthm,multirow,color,xcolor,cancel,bbm,tikz,mathrsfs,hyperref,authblk,stmaryrd,dsfont}
\usepackage[textheight=24cm,textwidth=16cm]{geometry}
\usetikzlibrary{decorations.markings}
\usetikzlibrary{shapes.multipart}

\hypersetup{
    colorlinks=true,
    linkcolor=blue,
    citecolor=cyan,
    filecolor=magenta,      
    urlcolor=green,
    pdftitle={Overleaf Example},
    pdfpagemode=FullScreen,
}

\newtheorem{lemma}{Lemma\setcounter{claimcounter}{0}}
\newtheorem{theorem}{Theorem\setcounter{claimcounter}{0}}

\newcounter{claimcounter}
\newtheorem*{claim*}{{\it Claim}}

\newcommand{\EM}[1]{{\it\textcolor{Maroon}{#1}}}
\newcommand{\EMM}[1]{{\textcolor{Maroon}{#1}}}

\def\lcm{\mathrm{lcm}}
\def\N{\mathbb{N}}
\def\F{\mathbb{F}}
\def\Z{\mathbb{Z}}
\def\X{\mathcal{X}}
\def\K{\mathcal{K}}
\def\D{\mathcal{D}}
\def\R{\mathcal{R}}

\def\S{\mathcal{S}}

\begin{document}

\title{There is no prime functional digraph: Seifert's proof revisited}

\author{Adrien Richard}

\maketitle

\begin{abstract}
A functional digraph is a finite digraph in which each vertex has a unique out-neighbor. Considered up to isomorphism and endowed with the directed sum and product, functional digraphs form a semigroup that has recently attracted significant attention, particularly regarding its multiplicative structure. In this context, a functional digraph $X$ divides a functional digraph $A$ if there exists a functional digraph $Y$ such that $XY$ is isomorphic to $A$. The digraph $X$ is said to be prime if it is not the identity for the product, and if, for all functional digraphs $A$ and $B$, the fact that $X$ divides $AB$ implies that $X$ divides $A$ or $B$. In 2020, Antonio E. Porreca asked whether prime functional digraphs exist, and in 2023, his work led him to conjecture that they do not. However, in 2024, Barbora Hudcová discovered that this result had already been proven by Ralph Seifert in 1971, in a somewhat forgotten paper. The terminology in that work differs significantly from that used in recent studies, the framework is more general, and the non-existence of prime functional digraphs appears only as a part of broader results, relying on (overly) technical lemmas developed within this general setting. The aim of this note is to present a much more accessible version of Seifert's proof — that no prime functional digraph exists — by using the current language and simplifying each step as much as possible.
\end{abstract}

\section{Introduction}

A deterministic, finite, discrete-time dynamical system is a function from a finite set (of states or configurations) to itself. Equivalently, it can be represented as a \EM{functional digraph}, that is, a finite directed graph $X$ with vertex set denoted \EM{$V(X)$}, in which each vertex $x$ has a unique out-neighbor, denoted \EM{$X(x)$}. We denote by \EM{$\F$} the set of all functional digraphs. In addition to being ubiquitous objects in discrete mathematics, such systems have numerous real-world applications: they are commonly used to model gene networks \cite{K69,T73,TK01,J02}, neural networks \cite{MP43,H82,G85}, reaction systems \cite{ehrenfeucht2007reaction}, social interactions \cite{PS83, GT83}, and more \cite{TA90,GM90}. 

\medskip
In this note, we consider functional digraphs up to isomorphism. This choice is motivated by the fact that many dynamical properties of interest are invariant under isomorphism: the number of fixed points, periodic points, limit cycles, the lengths of limit cycles, and so on. See, for example, \cite{FO89} in the context of random functional digraphs, and \cite{R19,G20} in the context of automata networks. An isomorphism class then corresponds to an unlabelled functional digraph, and for any $A, B \in \F$, we write \EM{$A = B$} to indicate that $A$ is isomorphic to $B$.

\medskip
There are two natural algebraic operations to construct larger systems from smaller ones. Given two functional digraphs $A$ and $B$, the \EM{addition} \EM{$A + B$} is defined as the disjoint union of $A$ and $B$, while the \EM{multiplication} \EM{$A \cdot B$} (or simply \EM{$AB$}) is the \EM{direct product} of $A$ and $B$: the vertex set is given by $V(AB) = V(A) \times V(B)$ and $(AB)(a,b) = (A(a), B(b))$ for all $(a,b) \in V(AB)$. Hence, $AB$ describes the parallel evolution of the dynamics represented by $A$ and $B$. Endowed with these two operations and taken up to isomorphism, the set of functional digraphs forms a commutative semiring called the \EM{semiring of functional digraphs}. The identity element for addition is the empty functional digraph, while for multiplication it is the cycle of length one, denoted \EM{$C_1$}. More generally, the cycle of length $\ell$ is denoted \EM{$C_\ell$}, and the sum of cycles consisting of $n$ cycles of length $\ell$ is denoted \EM{$nC_\ell$}. Then, for all $a,b \geq 1$, one can easily verify (see e.g., \cite{dennunzio2018polynomial}) that the product of $C_a$ and $C_b$ satisfies 
\begin{equation}\label{eq:cycle_product}
C_aC_b=\left(\frac{ab}{\lcm(a,b)}\right)C_{\lcm(a,b)}=\gcd(a,b)C_{\lcm(a,b)}.
\end{equation}

\medskip
The semiring of functional digraphs contains an isomorphic copy of $\N$ (by identifying each integer $n$ with $nC_1$) and, essentially because the product of cycles described above involves the least common multiple operation, it exhibits interesting properties that differ from those of $\mathbb{N}$ (or the semiring of polynomials). These properties have been the subject of recent studies \cite{dennunzio2018polynomial, dennunzio2019solving, GMP20, riva2022factorisation, DFMR23, dore2024roots, dore2024decomposition, naquin2024factorisation, dennunzio2024note, porreca2025injectivity, porreca2025solving, BCPR26}, leading to challenging open algebraic and complexity problems.

\medskip
To emphasize such differences, we first need some definitions. Let us say that $X \in \F$ is \EM{irreducible} if, for all $A,B \in \F$, $X = AB$ implies $A = C_1$ or $B = C_1$. Then there is no unique factorization into irreducibles: $2C_2 = C_2 C_2 = (2C_1) C_2$, while $C_2$ and $2C_1$ are irreducible. Now, let us say that $A$ \EM{divides} $B$, denoted \EM{$A \mid B$}, if $AX = B$ for some $X \in \F$; note that $X$ is not necessarily unique: the previous example shows that $C_2 X = 2C_2$ for $X = C_2$ or $X = 2C_1$. Finally, let us say that $X \in \F$ is \EM{prime} if 
\[
X \neq C_1 \textrm{ and, } \forall A,B \in \F, \quad X \mid AB \implies X \mid A \textrm{ or } X \mid B.
\]
In $\mathbb{N}$, irreducibility and primeness are equivalent, but this is far from true for $\F$: almost all functional digraphs are irreducible \cite{Dorigatti2017}, while Ralph Seifert proved in 1971 the following.

\begin{theorem}[\cite{seifert1971prime}]\label{thm:Seifert}
There is no prime functional digraph. 
\end{theorem}

Without knowledge of Seifert's paper, Antonio E. Porreca asked in 2020 whether prime functional digraphs exist. Johan Couturier and Marius Rolland completed their Master's theses, under Porreca's supervision, on this interesting and apparently new open problem. Subsequently, Porreca conjectured that no prime functional digraphs exist. While attending the AUTOMATA 2024 conference, Barbora Hudcov\'a became aware of this problem, discovered Seifert's paper, and communicated it to us. It must be said, in defence of Porreca, that this paper appears to have received very little attention. To our knowledge, prior to Hudcov\'a's communication, it had only three citations, all in algebraic contexts far from the recent studies on functional digraphs. Furthermore, the language used by Seifert differs significantly from that used above to state his result (for instance, functional digraphs are called unitary algebras).

\medskip
Theorem \ref{thm:Seifert} is actually only a small part of the results presented in Seifert's paper; however, it is the only result formally stated in the introduction and appears as one of the most important. Seifert works within a much broader framework, and by reading the paper, it is difficult to discern which arguments are truly necessary for the proof of Theorem \ref{thm:Seifert}. The aim of this note is therefore to provide a much more transparent proof of Theorem \ref{thm:Seifert} by simplifying each step as much as possible and by using the language of the recent studies cited above.

\medskip
This simplified proof is given in Section \ref{sec:no_prime}. It consists of three steps. The first shows that every prime functional digraph is connected, and the second shows that every prime functional digraph contains $C_1$. These first two steps are straightforward to understand and were independently obtained in Johan Couturier's Master's thesis, using rather different arguments; however, here we choose to follow Seifert's approach. The third, and most difficult step, is to show that every connected functional digraph containing $C_1$ is not prime. Here we use exactly the constructive proof of Seifert, from which it is difficult to draw intuition, and which for this reason appears as a deep inspiration. Throughout the proof, we develop each argument (which is not really the case in Seifert's paper) to make the proof as accessible as possible to a non-expert reader.

\medskip
In Section \ref{sec:Seifert_work}, we explain why the proof of Theorem \ref{thm:Seifert} is not very accessible in Seifert's paper. In brief, he works within a larger framework, and all the lemmas leading to the proof of Theorem \ref{thm:Seifert} are made as general as possible. Many of these lemmas are interesting in their own right but require technical definitions and arguments that are unnecessary for Theorem \ref{thm:Seifert}. We will present the statements of all these intermediate lemmas in the current language, both to justify our simplification effort and to highlight the profound results of this forgotten paper.

\paragraph{Preliminaries} 
Let $A,B\in\F$. We denote by \EM{$|A|$} the number of vertices in $A$. Obviously, for all $A,B\in\F$, we have $|AB|=|A||B|$. As mentioned above, \EM{$A=B$} means that $A$ and $B$ are isomorphic; formally, there exists a bijection $\phi:V(A)\to V(B)$ such that $\phi(A(a))=B(\phi(a))$ for all $a\in V(A)$, and $\phi$ is then called an \EM{isomorphism} from $A$ to $B$. We write \EM{$A\subseteq B$} to mean that $A$ is isomorphic to a subgraph of $B$. If $n$ is a positive integer, then \EM{$nA$} is the functional digraph consisting of the disjoint union of $n$ copies of $A$. A \EM{sum of cycles} is a functional digraph in which each connected component is a cycle. For positive integers $a,b$, we set \EM{$a \lor b = \lcm(a,b)$} and \EM{$a \land b = \gcd(a,b) = ab/(a \lor b)$}. The formula \eqref{eq:cycle_product}, that $C_aC_b = (a \land b) C_{a \lor b}$, will be used many times without explicit reference.

\section{Proof of Theorem \ref{thm:Seifert}}\label{sec:no_prime}

\subsection{Prime functional digraphs are connected}

We begin with an easy lemma. 

\begin{lemma}\label{lem:irr} 
Let $A,B$ be non-empty sums of cycles. {\em (a)}  $C_\ell\subseteq AB$ iff there exists $C_a\subseteq A$ and $C_b\subseteq B$ such that $a\lor b=\ell$. {\em (b)} If $a,b$ is the maximum length of a cycle in $A,B$, respectively, then each cycle of $AB$ has length at most $ab$. {\em (c)} If $A=C_\ell$ then each cycle of $AB$ has length at least $\ell$. {\em (d)} If $\ell$ is a prime power, then $C_\ell$ is irreducible.
\end{lemma} 

\begin{proof}
Since $A=\sum_{i=1}^n C_{a_i}$ and $B=\sum_{i=1}^m C_{b_i}$ for some positive integers $n$,$m$,$a_1$,$\dots$,$a_n$,$b_1$,$\dots$,$b_m$, by distributivity, we have
\begin{equation}\label{eq:prod_sum_cycles}
AB=\sum_{i=1}^n\sum_{j=1}^m C_{a_i}C_{b_j}=\sum_{i=1}^n\sum_{j=1}^m (a_i\land b_j)C_{a_i\lor b_j}.
\end{equation}
We deduce that $C_\ell\subseteq AB$ if and only if there are $a_i,a_j$ such that $a_i\lor b_j=\ell$. This proves (a). Let $a,b$ be the maximum length of a cycle in $A,B$, respectively. Then for all $1\leq i\leq n$ and $1\leq j\leq m$, we have $a_i\lor b_j\leq a_ib_i\leq ab$. This proves (b). If $A=C_\ell$ then, by \eqref{eq:prod_sum_cycles}, each cycle in $AB$ has length $\ell\lor b_j\geq \ell$ for some $1\leq j\leq m$. This proves (c). Let $p^\alpha$ be a prime power. Suppose that  $AB=C_{p^\alpha}$. By \eqref{eq:prod_sum_cycles} we have $n=m=1$, and $a_1\lor b_1=p^\alpha$. So $a_1$ and $b_1$ divide $p^\alpha$, and at least one of them is $p^\alpha$, say $a_1$. Then $|A|=p^\alpha$ and since $|A||B|=|C_{p^\alpha}|={p^\alpha}$, we have $|B|=1$, that is, $B=C_1$. Thus $C_{p^\alpha}$ is irreducible. This proves (d).  
\end{proof}

Let $A,B,X\in \F$. We denote by \EM{$[X]$} the \EM{cyclic part} of $X$, that is, the sum of the cycles contained in $X$; the number of cycles in $[X]$ is the number of connected components in $X$. A vertex  $(a,b)$ of $AB$ is in $[AB]$ if and only if $a$ is in a cycle of $A$ and $b$ is in a cycle of $B$, thus
\begin{eqnarray}%
[A\cdot B]&=&[A]\cdot [B]\label{eq:cyclic_part_1}\\
A\mid B&\Rightarrow& [A]\mid [B]\label{eq:cyclic_part_2}
\end{eqnarray}

\medskip
With this basic material, we can now prove the following. 

\begin{lemma}\label{lem:connected}
Every prime functional digraph is connected. 
\end{lemma}

\begin{proof}
Suppose that $X\in\F$ is disconnected, and let us prove that $X$ is not prime. Let $\ell$ be the maximum length of a cycle in $X$, and let $X_1$ be a connected component of $X$ containing a cycle of length $\ell$. Then $X=X_1+X_2$ for some non-empty $X_2\in\F$. Let $p$ be a prime strictly greater than $\ell$. Since $C_pC_p=pC_p=(pC_1)C_p$, we have 
\[
X\cdot pC_p= X_1C_pC_p+X_2(pC_1)C_p = C_p\cdot (C_pX_1+pX_2).
\]
We have $X\nmid C_p$ since $C_p$ is irreducible (Lemma \ref{lem:irr}(d)) and $X\neq C_p$ (since $X$ is not connected). Suppose, for a contradiction, that $X\mid  C_pX_1+pX_2$. Then, by \eqref{eq:cyclic_part_2}, $[X]\mid C_pC_\ell+p[X_2]$. Since $\ell$ and $p$ are coprime, $C_\ell C_p=C_{\ell p}$. Hence there exists a sum of cycles $Y$ such that 
\begin{equation}\label{eq:connected}
[X] Y=C_{p\ell}+p[X_2].
\end{equation}
We have $|[X]|\cdot |Y|=p(\ell+|[X_2]|)=p\cdot |[X]|$, thus $|Y|=p$. Suppose that $Y\neq C_p$, and let $q$ be the maximum length of a cycle in $Y$; hence $q<p$. By Lemma \ref{lem:irr}(b), each cycle of $[X] Y$ has length at most $\ell q<\ell p$, and this contradicts \eqref{eq:connected}. Thus $Y=C_p$. But then, by Lemma \ref{lem:irr}(c), each cycle of $[X]Y$ has length at least $p$, and since each cycle of $p[X_2]$ has length at most $\ell<p$, this again contradicts~\eqref{eq:connected}. Thus $X\nmid  C_pX_1+pX_2$, and this proves that $X$ is not prime.
\end{proof}

\subsection{Prime functional digraphs contain a cycle of length one}

We denote by \EM{$\F_\ell$} the set of connected functional digraphs with a cycle of length $\ell$. Hence $[X]=C_\ell$ iff $X\in \F_{\ell}$. Given $X\in \F_\ell$ and a positive integer $n$, the product $C_nX$ is completely described in \cite[Lemma 13]{naquin2024factorisation}, but the statement needs the introduction of several definitions. We choose to prove something weaker on $C_nX$, sufficient for our purpose, without any additional definition. The proof is straightforward, and it is sufficient to analyse a few examples to get it; if it appears lengthy, this is because we develop each argument.   

\begin{lemma}\label{lem:C_nX}
Let $n,\ell\geq 1$ and $X\in \F_\ell$. {\em (a)} We have $C_nX=(n\land \ell) X'$ for some $X'\in \F_{n\lor \ell}$. {\em (b)} If $n\mid \ell$ then $C_nX=n X$. 
\end{lemma}

\begin{proof}
If $n=1$, there is nothing to prove, so suppose that $n>1$. By \eqref{eq:cyclic_part_1} we have $[C_nX]=C_n[X]=C_nC_\ell=(n\land \ell)C_{n\lor \ell}$. We deduce that  $C_nX$ has $(n\land \ell)$ connected components, each in $\F_{n\lor \ell}$. Say, without loss of generality, that the vertices of $C_n$ are $0,1,\dots,n-1$ in order. In the remainder of this proof, we denote by $\oplus$ the sum modulo $n$, so $C_n(i)=i\oplus 1$ for all $0\leq i<n$, and we fix a vertex $x_0$ in $[X]$.  

\medskip
Let $Y_1$ be a connected component of $C_nX$, and let $V_1=V(Y_1)$. We will prove that there exists $0\leq i<n$ such that $(i,x_0)\in V_1$. Let any $(j,x)\in V_1$, and let $d$ be the length of the path from $x$ to $x_0$ in $X$, which exists since $X$ is connected and $x_0$ is in the cycle of $X$. Then $(j\oplus d,x_0)\in V_1$, so it suffices to take $i=j\oplus d$.  

\medskip
Suppose now that $C_nX$ has a connected component $Y_2$, distinct from $Y_1$, and let $V_2=V(Y_2)$. To prove the first assertion of the lemma, it is sufficient to prove that $Y_1=Y_2$. We have $(i,x_0)\in V_1$ and the argument above shows that $(i\oplus k,x_0)\in V_2$ for some $1\leq k<n$. We then define the map $\phi$ on $V_1$ as follows: for all $(j,x)\in V_1$, 
\[
\phi(j,x)=(j\oplus k,x).
\]

\medskip
Let us prove that $\phi$ is an isomorphism from $Y_1$ to $Y_2$. We first prove that $\phi(V_1)\subseteq V_2$. For that, we take any $(j,x)\in V_1$ and we prove that $(j\oplus k,x)\in V_2$ by induction on the length $d(j,x)$ of the path from $(j,x)$ to $(i,x_0)$ in $Y_1$, which exists since $Y_1$ is connected and $(i,x_0)$ belongs to the cycle of $Y_1$. If $d(j,x)=0$, then $(j,x)=(i,x_0)$ and $(i\oplus k,x_0)\in V_2$ by the choice of $k$. Suppose that $d(j,x)>0$. Then $Y_1(j,x)=(j\oplus 1,X(x))$ and since $d(Y_1(j,x))=d(j,x)-1$, by induction, $(j\oplus 1\oplus k,X(x))\in V_2$. Since $(C_n X)(j\oplus k,x)=(j\oplus 1\oplus k,X(x))$ we have $(j\oplus k,x)\in V_2$, as desired. So $\phi(V_1)\subseteq V_2$ and since $\phi$ is an injection, we have $|V_1|\leq |V_2|$. By switching the roles of $Y_1$ and $Y_2$ and using the same argument, we get that $|V_2|\leq |V_1|$. Thus $\phi$ is a bijection from $V_1$ to $V_2$. Furthermore, it is an isomorphism from $Y_1$ to $Y_2$ since, for all $(j,x)\in V_1$, we have  
\[
\phi(Y_1(j,x))=\phi(j\oplus 1,X(x))=(j\oplus k\oplus 1,X(x))=(C_nX)(j\oplus k,x)=(C_nX)(\phi(j,x))=Y_2(\phi(j,x)).
\]

\medskip
Suppose now that $n\mid \ell$. Since $n\land \ell=n$, we have already proven that $C_n X=n Y_1$. It remains to prove that $Y_1=X$. Let $\phi$ be the map from $V_1$ to $V(X)$ defined as follows: for all $(j,x)\in V_1$, $\phi(j,x)=x$. Let us prove that $\phi$ is a bijection. Since $|Y_1|=|X|$ it is sufficient to prove that $\phi$ is surjective. For that, we take any $x\in V(X)$ and prove that $(j,x)\in V_1$ for some $0\leq j<n$ by induction on the length $d(x)$ of the path from $x$ to $x_0$ in $X$, which exists since $X$ is connected and $x_0$ belongs to the cycle of $X$. If $d(x)=0$ then $x=x_0$ and $(i,x_0)\in V_1$. Suppose that $d(x)>0$. Since $d(X(x))=d(x)-1$, by induction, there exists $0\leq j<n$ such that $(j\oplus 1,X(x))\in V_1$, and since  $(C_nX)(j,x)=(j\oplus 1,X(x))$ we have $(j,x)\in V_1$, as desired. Thus $\phi$ is a bijection from $V_1$ to $V(X)$, and it is an isomorphism from $Y_1$ to $X$ since, for all $(j,x)\in V_1$, we have $\phi(Y_1(j,x))=\phi(j\oplus 1,X(x))=X(x)=X(\phi(j,x))$.
\end{proof}

We now have a very simple argument to prove the following, which, except for the use of Lemma~\ref{lem:C_nX}, only relies on the cyclic part. 

\begin{lemma}\label{lem:F_1}
Every connected prime functional digraph has a cycle of length 1. 
\end{lemma}

\begin{proof}
Let $\ell>1$ and $X\in\F_\ell$. Let $p$ be a prime which divides $\ell$, and let $\alpha$ be the largest integer such that $p^{\alpha}\mid \ell$. If $X=C_{p^\alpha}$ then 
\[
C_{p^\alpha}\cdot p C_{p^{\alpha+1}}=p(C_{p^\alpha}C_{p^{\alpha+1}})=p^{\alpha+1}C_{p^{\alpha+1}}= C_{p^{\alpha+1}}\cdot C_{p^{\alpha+1}}, 
\]
thus $X$ is not prime since $C_{p^{\alpha+1}}$ is irreducible (Lemma \ref{lem:irr}(d)) and $X\neq C_{p^{\alpha+1}}$. So suppose that $X\neq C_{p^\alpha}$. By Lemma \ref{lem:C_nX}(a), there exists $X'\in\F_{\ell^2}$, such that $XC_{\ell^2}=\ell X'$, and by Lemma \ref{lem:C_nX}(b) we have $\ell X'=C_\ell X'$. Since $C_\ell=C_{p^\alpha}\cdot C_{\ell/p^{\alpha}}$, we obtain 
\[
XC_{\ell^2}=C_{p^\alpha}\cdot (C_{\ell/p^{\alpha}} X').
\] 
Then $X\nmid C_{p^\alpha}$ since $C_{p^\alpha}$ is irreducible (Lemma \ref{lem:irr}(d)) and $X\neq C_{p^\alpha}$. Suppose that $X\mid C_{\ell/p^{\alpha}} X'$. Then, by \eqref{eq:cyclic_part_2}, we have $C_\ell\mid C_{\ell/p^{\alpha}} C_{\ell^2}$. Hence there exists a sum of cycles $Y$ such that $C_\ell Y= C_{\ell/p^{\alpha}} C_{\ell^2}=(\ell/p^{\alpha})C_{\ell^2}$. We deduce from Lemma \ref{lem:irr}(a) that $C_{\ell^2}\subseteq Y$. Thus 
\[
|C_\ell Y|\geq \ell^3> \frac{\ell^3}{p^{\alpha}}=|C_{\ell/p^{\alpha}} C_{\ell^2}|,
\]
a contradiction. Thus $X$ is not prime.  
\end{proof}

\subsection{Seifert's proof that $\F_1$ has no prime functional digraphs}

By Lemmas~\ref{lem:connected} and \ref{lem:F_1}, every prime functional digraph is in $\F_1$. We now present Seifert's ingenious constructive proof that $\F_1$ has no prime, in Lemma~\ref{lem:F_1_no_prime} below (there is just an insignificant change in the definition of $B$). This concludes the proof that there is no prime functional digraph (Theorem~\ref{thm:Seifert}). If the proof that $\F_1$ has no prime appears longer than that of Seifert's paper, this is because we develop each argument. 

\medskip
We first need some preliminaries. Regarding $X\in\F$ as a function from $V(X)$ to itself, we denote by \EM{$X^k$} the $k$-fold composition of $X$ with itself: $X^0(x)=x$, and $X^k(x)=X(X^{k-1}(x))$ for all $k\geq 1$. Let $X\in \F_1$, and let $\chi$ be the vertex in the cycle of length $1$ in $X$. In other words, $\chi$ is the unique fixed point of $X$. For all $x\in V(X)$, we denote by \EM{$d_X(x)$} the length of the path from $x$ to $\chi$ (which exists since $X$ is connected). Equivalently, $d_X(x)$ is the smallest integer $d$ such that $X^d(x)=\chi$; in particular, $d_X(\chi)=0$ and if $x\neq\chi$ then $d_X(X(x))=d_X(x)-1$. For $d\geq 0$, we denote by \EM{$X_{\mid d}$} the subdigraph of $X$ induced by the set of vertices $x\in V(X)$ with $d_X(x)\leq d$; we have $X_{\mid d}\in\F_1$. The \EM{height} of $X$, denoted \EM{$d(X)$}, is the maximum of $d_X(x)$ for $x\in V(X)$. Let $Y\in\F_1$ and $(x,u)\in V(XY)$. We easily check that $d_{XY}(x,y)=\max(d_X(X),d_Y(Y))$ and thus 
\begin{equation}\label{eq:max}
d(XY)=\max(d(X),d(Y)).
\end{equation}

\begin{lemma}\label{lem:F_1_no_prime}
There is no prime functional digraph in $\F_1$. 
\end{lemma}

\begin{proof}[{\bf Proof}]
Let $X\in \F_1$ with $X\neq C_1$. We will construct $A,B,Y\in \F$ such that $XY=AB$ but $X\nmid A$ and $X\nmid B$, thus proving that $X$ is not prime. The equality $XY=AB$ will be obtained with an explicit isomorphism $\phi$ from $XY$ to $AB$. Let $\chi$ be the fixed point of $X$ and let $d=d(X)$; we have $d\geq 1$ since $X\neq C_1$. 
\smallskip
\begin{quote}
{\bf Definition of \boldmath$A$\unboldmath.} Let $P$ be the functional digraph defined by:
\begin{itemize}
\item $V(P)=\{0,1,\dots,d\}$,
\item $P(i)=\max(i-1,0)$ for all $i\in V(P)$. 
\end{itemize}
Hence $P\in\F_1$, its fixed point is $0$, and it is isomorphic to the functional digraph obtained from a path of length $d$ by adding a loop on the last vertex. Then $A$ is the functional digraph obtained from $XP$ by adding a new vertex $u$, not in $XP$, and an edge from $u$ to $(\chi,d-1)$:
\begin{itemize}
\item $V(A)=V(XP)\cup\{u\}$,
\item $A(x,i)=(X(x),P(i))$, for all $(x,i)\in V(XP)$,
\item $A(u)=(\chi,d-1)$.
\end{itemize}
An illustration is given in Figure \ref{fig:A}.
\end{quote}

\input{figA.tex}
\input{figB.tex}
\input{figY.tex}

\medskip
\noindent
{\it Proof that $X\nmid A$.} Indeed, $|A|/|X|=|P|+(1/|X|)$ is not an integer since $|X|>1$. 

\smallskip
\begin{quote} {\bf Definition of \boldmath$B$\unboldmath.} Let $t$ be a new vertex, not in $X$. Let $\X$ be the functional digraph obtained from $X$ by adding $t$ and an edge from $t$ to an arbitrary vertex $\hat x\in V(X)$ with $d_X(\hat x)=d$; hence $\X\in \F_1$, its fixed point is $\chi$, and the important fact is that $d(\X)=d_{\X}(t)=d+1$. So we only need the following specifications for $\X$:
\begin{itemize}
\item $V(\X)=V(X)\cup\{t\}$,
\item $\X(x)=X(x)$ for all $x\in V(X)$,
\item $d_{\X}(t)=d+1$.
\end{itemize}
Then $B$ is defined as follows: 
\begin{itemize}
\item $V(B)=V_1\times\cdots\times V_{d+1}$, where $V_i=V(\X_{\mid i})$ for $1\leq i\leq d+1$,
\item $B(b)_i=\X(b_{i+1})$, $1\leq i\leq d$, and $B(b)_{d+1}=t$, for all $b=(b_1,\dots,b_{d+1})\in V(B)$. 
\end{itemize}
An illustration is given in Figure \ref{fig:B}. Let us make some remarks. Firstly, $B$ is well defined since, for all $0\leq i\leq d$, we have $b_{i+1}\in V_{i+1}$ and thus $\X(b_{i+1})\in V_i$. Secondly, $B(b)_i=X(b_{i+1})$ for $1\leq i<d$; this will be used many times. Thirdly, $B$ does not depend on component $1$: if $b,b'\in V(B)$ only differ in $b_1\neq b'_1$ then $B(b)=B(b')$. Lastly, for all $1\leq i,j\leq d+1$ with $j\neq i+1$, $B_i$ does not depend on $j$: if $b,b'\in V(B)$ only differ in $b_j\neq b'_j$ then $B(b)_i=B(b')_i$. The last two remarks will be used to prove that $X\mid AB$. 
\end{quote}

\medskip
\noindent
{\it Proof that $X\nmid B$.} We first need some properties on $B$. Let $\beta\in V(B)$ be defined by 
\[
\beta_i=\X^{d-i+1}(t)\textrm{ for }1\leq i\leq d+1. 
\]
We have $B(\beta)_{d+1}=t=\X^0(t)=\beta_{d+1}$, and for $1\leq i< d+1$, 
\begin{equation}\label{eq:i-1}
B(\beta)_i=\X(\beta_{i+1})=\X(\X^{d-i}(t))=\X^{d-i+1}(t)=\beta_i. 
\end{equation}
Thus $\beta$ is a fixed point of $B$. Let us prove that
\begin{equation}\label{eq:k}
B^k(b)_i=\beta_i\textrm{ for all }b\in V(B)\textrm{, }1\leq i\leq d+1\textrm{ and }k\geq d-i+2. 
\end{equation}
We proceed by induction on $i$ from $d+1$ to $1$. We have $B^k(b)_{d+1}=t=\beta_{d+1}$ for all $k\geq 1$. Suppose now that \eqref{eq:k} holds for $1<i\leq d+1$, and let $k\geq d-i+2$. Then $B^{k+1}(b)_{i-1}=B(B^k(b))_{i-1}=\X(B^k(b)_i)=\X(\beta_i)=\beta_{i-1}$, where the before last equality is obtained using the induction hypothesis, and the last equality is given by \eqref{eq:i-1}. This completes the induction step. We deduce that $B^{d+1}(b)=\beta$, hence $B\in\F_1$ and $d(B)\leq d+1$. We now prove the following equivalence, which shows that $d(B)=d+1$ by characterising the vertices with height at most~$d$:
\begin{equation}\label{eq:d}
B^d(b)=\beta\iff b_{d+1}=t.
\end{equation}
Indeed, we easily check, by induction on $k$ from $1$ to $d$, that $B^k(b)_1=\X^k(b_{k+1})$. In particular, $B^d(b)_1=\X^d(b_{d+1})$ and, by \eqref{eq:k}, we have $B^d(b)_i=\beta_i$ for all $1<i\leq d+1$. Thus $B^d(b)=\beta$ iff $B^d(b)_1=\beta_1$, that is, $\X^{d}(b_{d+1})=\X^d(t)$. So $B^d(b)=\beta$ if $b_{d+1}=t$. If $b_{d+1}\neq t$ then $\X^d(b_{d+1})=X^d(b_{d+1})=\chi\neq \X^d(t)$, and this proves \eqref{eq:d}. Furthermore, we have $B^{d-1}(b)_1=\X^{d-1}(b_d)=X^{d-1}(b_d)$ and, using very similar arguments, we get that if $b_{d+1}=t$, then $B^{d-1}(b)_i=\beta_i$ for $1<i\leq d+1$. We deduce that, for all $b,b'\in V(B)$, 
\begin{equation}\label{eq:d-1}
b_{d+1}=b'_{d+1}=t\quad\Longrightarrow\quad \left(B^{d-1}(b)=B^{d-1}(b')~\iff~X^{d-1}(b_d)=X^{d-1}(b'_d)\right).
\end{equation}

\medskip
Suppose now, for a contradiction, that $X\mid B$, and let $Z\in \F$ such that there is an isomorphism, say $\psi$, from $XZ$ to $B$. Then $C_1[Z]=[X][Z]=[XZ]=[B]=C_1$, thus $[Z]=C_1$, that is, $Z\in \F_1$. For $E\in \F_1$, let $\sim_E$ be the equivalence relation on $V(E_{\mid d})$ defined by 
\[
e\sim_E e' ~ \iff ~ E^{d-1}(e)=E^{d-1}(e'),
\]
and let $n_E$ be the number of equivalent classes of $\sim_E$. By \eqref{eq:d}, $V(B_{\mid d})$ is the set of $b\in V(B)$ with $b_{d+1}=t$, and since $X=X_{\mid d}$, we have $V(X_{\mid d})=V_d$. Consequently, by $\eqref{eq:d-1}$, we have 
\[
n_X=n_B,
\]
which is the key observation. Now, for all $(x,z),(x',z')\in V(X)\times V(Z_{\mid d})=V((XZ)_{\mid d})$, we have $(x,z)\sim_{XZ} (x',z')$ iff $x\sim_X x'$ and $z\sim_Z z'$. Thus $n_{XZ}=n_X\cdot n_Z$. Let $k\geq 1$ and $(x,z)\in V(XZ)$. We easily prove, by induction on $k$, that $\psi((XZ)^k(x,z))=B^k(\psi(x,z))$. We deduce that $(x,z)\sim_{XZ} (x',z')$ iff $\psi(x,z)\sim_B \psi(x',z')$, and thus $n_{XZ}=n_B$. So $n_X\cdot n_Z=n_B$ and since $n_X=n_B$, we have $n_Z=1$. Let $\zeta$ be the fixed point of $Z$. That $n_Z=1$ means that $Z^{d-1}(z)=Z^{d-1}(\zeta)=\zeta$ for all $z\in V(Z_{\mid d})$. Hence $d(Z_{\mid d})\leq d-1$, thus $d(Z)\leq d-1$. We deduce from \eqref{eq:max} that $d(B)=\max(d(X),d(Z))=d$, and this contradicts \eqref{eq:d}. This proves that $X\nmid B$.

\smallskip
\begin{quote}
{\bf Definition of \boldmath$Y$\unboldmath.}
Let $V_\chi(B)$ be the set of $b\in V(B)$ such that $b_d=\chi$. Then $Y$ is obtained from $PB$ by adding the vertices in $V_\chi(B)$ and an edge from $b$ to $(d-1,B(b))\in V(PB)$ for each $b\in V_\chi(B)$:
\begin{itemize}
\item $V(Y)=V(PB)\cup V_\chi(B)$, 
\item $Y(i,b)=(P(i),B(d))$ for all $(i,b)\in V(PB)$,
\item $Y(b)=(d-1,B(b))$ for all $b\in V_\chi(B)$. 
\end{itemize}
An illustration is given in Figure \ref{fig:Y}.
\end{quote}

\begin{quote}
{\bf Definition of \boldmath$\phi$\unboldmath.}
Given $b\in V(B)$ and $x\in V_i$, $1\leq i\leq d+1$, we denote by $b^{i,x}$ the vertex in $B$ obtained from $b$ by replacing its $i$th component by $x$: $(b^{i,x})_i=x$ and $(b^{i,x})_j=b_j$ for all $j\neq i$. Furthermore, we set $b^{0,x}=b$ and $b_0=\chi$. Then $\phi$ is the map on $V(XY)$ defined by: for all $x\in V(X)$, $(i,b)\in V(PB)$, and $b'\in V_\chi(B)$, 
\begin{itemize}
\item $\phi(x,(i,b))=((b_i,i),b^{i,x})$ if $d_X(x)\leq i$, 
\item $\phi(x,(i,b))=((x,i),b)$ if $d_X(x)> i$, 
\item $\phi(x,b')=(u,b'^{d,x})$. 
\end{itemize}
An illustration is given in Figure \ref{fig:XYAB}.
\end{quote}

\input{figXYAB.tex}

\medskip
\noindent
{\it Proof that $\phi$ is an isomorphism from $XY$ to $AB$.} This is a straightforward computation, which is not developed in Seifert's paper for that reason. We give the proof here for completeness. We easily check that each image of $\phi$ is in $V(AB)$, and that $\phi$ is injective. Furthermore, $\phi$ is a bijection since, $|V_\chi(B)|=|B|/|V_d|=|B|/|X|$ and thus 
\[
|X||Y|=|X|(|P||B|+|B|/|X|)=|X||P||B|+|B|=(|X||P|+1)|B|=|A||B|.
\] 

\medskip
Next, we need the following property: for all $b\in V(B)$, $i\in V(P)$ and $x\in V(X)$, 
\begin{equation}\label{eq:phi}
B(b)^{P(i),X(x)}=B(b^{i,x}).
\end{equation}
Indeed, if $i\leq 1$, then $P(i)=0$ and $B(b)^{0,X(x)}=B(b)=B(b^{0,x})=B(b^{1,x})$, where the last equality holds since $B$ does not depend on component $1$. Suppose now that $i\geq 2$. Then 
\[
\left(B(b)^{P(i),X(x)}\right)_{i-1}=\left(B(b)^{i-1,X(x)}\right)_{i-1}=X(x)=X((b^{i,x})_i)=B(b^{i,x})_{i-1}, 
\]
and for all $1\leq j\leq d+1$ with $j\neq i-1$, we have 
\[
\left(B(b)^{P(i),X(x)}\right)_j=\left(B(b)^{i-1,X(x)}\right)_j=B(b)_j=B(b^{i,x})_j,
\]
where the last equality holds since $B_j$ does not depend on component $i$. This proves \eqref{eq:phi}. 

\medskip
We can now prove $\phi$ is an isomorphism, that is, $\phi(XY(z))=AB(\phi(z))$ for all $z\in V(XY)$. Let $x\in V(X)$, $(i,b)\in V(PB)$, and $b'\in V_\chi(B)$. We consider four cases. 
\begin{itemize}
\item
Suppose that $d_X(x)\leq i$. Then $d_X(X(x))\leq P(i)$. Indeed, if $x=\chi$, then $d_X(X(x))=d_X(\chi)=0\leq P(i)$. Otherwise, $1\leq d_X(x)$, and thus $d_X(X(x))=d_X(x)-1\leq i-1=P(i)$. Consequently,
\[
\begin{array}{lll}
AB(\phi(x,(i,b)))&=AB((b_i,i),b^{i,x})&=((X(b_i),P(i)),B(b^{i,x}))\\
\phi(XY(x,(i,b)))&=\phi(X(x),(P(i),B(b)))&=((B(b)_{P(i)},P(i)),B(b)^{P(i),X(x)}).
\end{array}
\]
So, by \eqref{eq:phi}, we have to prove $B(b)_{P(i)}=X(b_i)$. If $i=0$, then $B(b)_{P(0)}=B(b)_0=\chi=X(b_0)$, and otherwise $B(b)_{P(i)}=B(b)_{i-1}=X(b_i)$, so we are done. 
\item
Suppose that $d_X(x)> i$ and $d_X(X(x))>P(i)$. Then we are done since 
\[
\begin{array}{lll}
AB(\phi(x,(i,b)))&=AB((x,i),b)&=((X(x),P(i)),B(b))\\
\phi(XY(x,(i,b)))&=\phi(X(x),(P(i),B(b)))&=((X(x),P(i)),B(b)).
\end{array}
\]
\item
Suppose that $d_X(x)> i$ and $d_X(X(x))\leq P(i)$. If $i\geq 1$, then $d_X(X(x))=d_X(x)-1>i-1=P(i)$, a contradiction, thus $i=0=P(i)$. So $X(x)=\chi$ and we are done since 
\[
\begin{array}{llll}
AB(\phi(x,(0,b)))&=AB((x,0),b)&=((X(x),0),B(b))&=((\chi,0),B(b))\\
\phi(XY(x,(0,b)))&=\phi(\chi,(0,B(b)))&=((B(b)_0,0),B(b)^{0,\chi})&=((\chi,0),B(b)).
\end{array}
\]
\item
Finally, since $d_X(X(x))\leq d-1$, we have 
\[
\begin{array}{lll}
AB(\phi(x,b'))&=AB(u,b'^{d,x})&=((\chi,d-1), B(b'^{d,x}))\\ 
\phi(XY(x,b'))&=\phi(X(x),(d-1,B(b')))&=((B(b')_{d-1},d-1),B(b')^{d-1,X(x)}).
\end{array}
\]
So we are done using \eqref{eq:phi} and the fact that $B(b')_{d-1}=\chi$ if $d=1$, and $B(b')_{d-1}=X(b'_d)=X(\chi)=\chi$ otherwise.
\end{itemize}
This proves the lemma.
\end{proof}


\section{Highlights on Seifert's original work}\label{sec:Seifert_work}

Seifert works on digraph families that are more general than functional digraphs. Given two digraphs $A,B$, $A=B$ still means that $A$ and $B$ are isomorphic, and we still denote $A\cdot B$ (or $AB$) as the direct product between $A$ and $B$: the vertex set of $AB$ is $V(A)\times V(B)$, and for all edges $(a,a')$ in $A$ and $(b,b')$ in $B$, $AB$ has an edge from $(a,b)$ to $(a',b')$. Then $A\mid B$ means that there exists a digraph $X$ such that $AX=B$. Seifert's definition of primeness is then not defined only for semirings of digraphs but more generally for any set $\K$ of digraphs (possibly including infinite digraphs), as follows: a digraph $X$ is \EM{$\K$-prime} if: (i) $X\in\K$, (ii) $X$ has more than one vertex, and (iii) for all $A,B\in\K$, if $X\mid AB$, then $X\mid A$ or $X\mid B$ \footnote{that is, if $XY=AB$ for some digraph $Y$, then $XY'=A$ for some digraph $Y'$ or $XY''=B$ for some digraph $Y''$, but $Y,Y'$ and $Y''$ are not necessarily in $\K$.}. The set of $\K$-prime digraphs is denoted \EM{$\Pr(\K)$}.

\medskip 
Seifert's paper has three chapters. The first chapter concerns infinite digraphs and is not connected to Theorem~\ref{thm:Seifert}, so we skip it. The other two chapters form the major part of the paper and treat four sets of finite digraphs: the set \EM{$\D$} of digraphs with minimum out-degree at least one, the set \EM{$\S$} of digraphs with minimum out- and in-degree at least one, the set \EM{$\R$} of digraphs with a loop on each vertex, and the set $\F$ of functional digraphs. So $\R \subseteq \S \subseteq \D$ and $\F \subseteq \D$. For $\K \in \{\R, \S, \F, \D\}$, let \EM{$\K_C$} be the set of connected digraphs in $\K$.

\medskip
The second chapter proves Theorem 2.7, which describes the number of connected components in $AB$ when $A,B\in \D_C$. For that, Seifert introduces a tricky digraph parameter as follows. Let $\EMM{A^1}=A$, and let $\EMM{A^{-1}}$ be the digraph obtained from $A$ by reversing the direction of every edge. Let \EM{$A\circ B$} be the digraph on $V(A)\cup V(B)$ with an edge from $a$ to $b$ if there exists $z$ such that $(a,z),(z,b)$ are edges of $A,B$, respectively. Given a sequence $\tau\in\{-1,1\}^*$ of length $k\geq 1$, let \EM{$w(\tau)$} be the sum of the elements of $\tau$, let \EM{$\tau_1$} be the first element of $\tau$; if $k>1$, we denote by \EM{$\tau'$} the sequence of length $k-1$ obtained from $\tau$ by deleting the first element. Then, \EM{$\tau(A)$} is the digraph inductively defined as follows: for $k=1$, $\tau(A)=A^{\tau_1}$, and for $k>1$, $\tau(A)= A^{\tau_1}\circ\tau'(A)$. Next, for $a\in V(A)$, let \EM{$K_A(a)$} be the set of $w(\tau)$ for non-empty $\tau\in\{-1,1\}^*$ such that $(a,a)$ is an edge of $\tau(A)$, and let \EM{$K(A)$} be the union of the sets $K_A(a)$ for $a\in V(A)$. Lemma 2.5 shows that $K_A(a)=K(A)$ for every $a\in V(A)$ and that $K(A)$ is a subgroup of $\Z$ distinct from $\{0\}$. Then the mentioned tricky parameter, denoted \EM{$k(A)$}, is the (unique) non-negative generator of $K(A)$, and Theorem 2.7 asserts that, given $A,B\in \D_C$, the number of connected components in $AB$ is $k(A)\land k(B)$. It takes some thought to digest the definition of $k(A)$ and to understand that, if $A$ is a connected functional digraph, then $k(A)$ is simply the length of the cycle in $A$. Then, if $A,B\in \F_C$, that $AB$ contains $k(A)\land k(B)$ connected components is an immediate consequence of the easy formulas \eqref{eq:cycle_product} and \eqref{eq:cyclic_part_1}. The proof of Theorem 2.7 relies on the fact that $K(A),K(B)$ are subgroups of $\Z$ and uses some arguments that are present in Lemma~\ref{lem:C_nX}. As a consequence of the proof, Seifert gets Corollary~2.8, that all the connected components of $C_nA$ are pairwise isomorphic, which thus generalises the first assertion in Lemma~\ref{lem:C_nX} to the class $\D_C$. But while Lemma~\ref{lem:C_nX} is an easy exercise, this generalisation is clearly non-trivial.

\medskip
The first goal of the third chapter is to prove that $\Pr(\F)=\Pr(\F_C)=\emptyset$. While only $\Pr(\F)=\emptyset$ (Theorem~\ref{thm:Seifert}) is mentioned in the introduction, the motivation for proving $\Pr(\F)=\Pr(\F_C)$ seems to be a similar equality previously obtained for the family of finite and infinite functional digraphs. Each step of the proof is as general as possible. Seifert first proves that, for $\K\in\{\S,\F,\D\}$ we have $\Pr(\K)\subseteq \K_C$, which implies $\Pr(\K)\subseteq \Pr(\K_C)$; this generalises Lemma~\ref{lem:connected} by including the classes $\S$ and $\D$. As in the proof of Lemma~\ref{lem:connected}, Seifert considers the product $C_pX$ for some large prime $p$ and the fact that $C_pC_p=pC_p=(pC_1)C_p$. But Lemma~\ref{lem:connected} is only based on basic observations concerning the cyclic part, while Seifert uses the rather complex material of the second chapter. The next result is that, for $\K\in\{\S,\F,\D\}$, if $X\in \Pr(\K)$ then $k(X)=1$; this generalises Lemma~\ref{lem:F_1} by including the classes $\S$ and $\D$. As in the proof of Lemma~~\ref{lem:F_1}, Seifert supposes, by contradiction, that $\ell=k(X)>1$, and considers the product $C_{\ell^2}X$. He then uses Corollary~2.8 (instead of Lemma~\ref{lem:C_nX}) and many additional technical arguments, while the proof of Lemma~\ref{lem:F_1} is very simple. Our simplification work is significant here. Then Seifert proves that $\Pr(\K_C)\subseteq \Pr(\K)$ for $\K\in\{\S,\R,\F,\D\}$, and thus $\Pr(\K_C)=\Pr(\K)$ for $\K\in\{\S,\F,\D\}$, which is not necessary to prove that $\Pr(\F)=\emptyset$. To conclude that $\Pr(\F)=\Pr(\F_C)=\emptyset$, he finally gives the constructive proof that $\Pr(\F_1)=\emptyset$, which is translated in our language in Lemma~\ref{lem:F_1_no_prime}. All of these occupy a major part of the chapter. The second and last goal is to prove that the $\R$-prime digraphs are exactly the irreducible members of $\R_C$ and the set of $pC_1$ for primes $p$.  

\paragraph{Acknowledgments} This work has been partially funded by the HORIZON-MSCA-2022-
SE-01 project 101131549 ‘‘ACANCOS’’ project, and the
ANR-24-CE48-7504 ‘‘ALARICE’’ project. I thank Antonio E. Porreca and Marius Rolland for their feedback. I also thank Eric Gibelin, Olivier Kreis, Maxime Réthoré, and Raphaël Soubeyran for their support.

\bibliographystyle{plain}
\bibliography{BIB}

\end{document}

%% file: figA.tex
\begin{figure}[p]
\[
\def\y{1}
\def\x{1}
\begin{array}{cccc}
\begin{array}{c}
\begin{tikzpicture}[every node/.style={outer sep=1,inner sep=0}]
\node (x) at (0,{1*\y}){\small{$\hat x $}};
\node (f) at (0,{0*\y}){\small{$\chi$}};
\path[decoration={markings,mark=at position 0.5 with {\arrow{>}}}]
(x) edge[postaction={decorate}] (f)
;
\draw[decoration={markings,mark=at position 0.55 with {\arrow{>}}},postaction={decorate}]
(f.-112) .. controls ({0-0.5},{0-0.7}) and ({{0+0.5}},{0-0.7}) .. (f.-68);
\end{tikzpicture}
\end{array}
&
\begin{array}{c}
\begin{tikzpicture}[every node/.style={outer sep=1,inner sep=0}]
\node (1) at (0,{1*\y}){\small{$1$}};
\node (0) at (0,{0*\y}){\small{$0$}};
\path[decoration={markings,mark=at position 0.5 with {\arrow{>}}}]
(1) edge[postaction={decorate}] (0)
;
\draw[decoration={markings,mark=at position 0.55 with {\arrow{>}}},postaction={decorate}]
(0.-112) .. controls ({0-0.5},{0-0.7}) and ({{0+0.5}},{0-0.7}) .. (0.-68);
\end{tikzpicture}
\end{array}
&
\begin{array}{c}
\begin{tikzpicture}[every node/.style={outer sep=1,inner sep=0}]
\node (b1) at ({-1*\x},{1*\y}){\small{$\hat x 1$}};
\node (f1) at ({+0*\x},{1*\y}){\small{$\chi 1$}};
\node (b0) at ({+1*\x},{1*\y}){\small{$\hat x  0$}};
\node (f0) at ({+0*\x},{0*\y}){\small{$\chi 0$}};
\path[decoration={markings,mark=at position 0.5 with {\arrow{>}}}]
(b1) edge[postaction={decorate}] (f0)
(f1) edge[postaction={decorate}] (f0)
(b0) edge[postaction={decorate}] (f0)
;
\draw[decoration={markings,mark=at position 0.55 with {\arrow{>}}},postaction={decorate}]
(f0.-112) .. controls ({0-0.5},{0-0.7}) and ({{0+0.5}},{0-0.7}) .. (f0.-68);
\end{tikzpicture}
\end{array}
&
\begin{array}{c}
\begin{tikzpicture}[every node/.style={outer sep=1,inner sep=0}]
\node (b1) at ({-1.5*\x},{1*\y}){\small{$\hat x 1$}};
\node (f1) at ({-0.5*\x},{1*\y}){\small{$\chi 1$}};
\node (b0) at ({+0.5*\x},{1*\y}){\small{$\hat x 0$}};
\node (u)  at ({+1.5*\x},{1*\y}){\small{\textcolor{blue}{$u$}}};
\node (f0) at ({+0.0*\x},{0*\y}){\small{$\chi 0$}};
\path[decoration={markings,mark=at position 0.5 with {\arrow{>}}}]
(b1) edge[postaction={decorate}] (f0)
(f1) edge[postaction={decorate}] (f0)
(b0) edge[postaction={decorate}] (f0)
(u)  edge[postaction={decorate},blue] (f0)
;
\draw[decoration={markings,mark=at position 0.55 with {\arrow{>}}},postaction={decorate}]
(f0.-112) .. controls ({0-0.5},{0-0.7}) and ({{0+0.5}},{0-0.7}) .. (f0.-68);
\end{tikzpicture}
\end{array}
\\
X&P&XP&A
\end{array}
\]
{\caption{\label{fig:A} Illustration of $A$, obtained from $XP$ by adding the blue vertex $(u$) and the blue edge.}}
\end{figure}
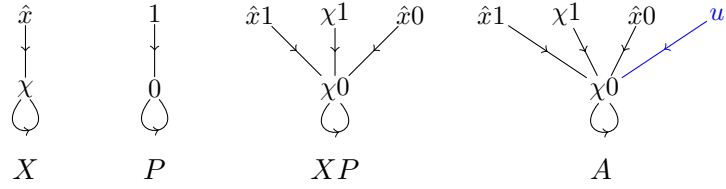

%% file: figB.tex
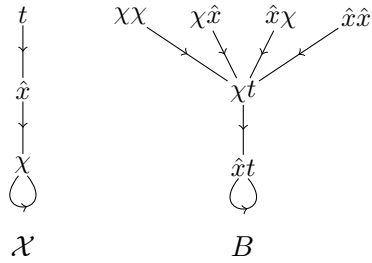
\begin{figure}[p]
\[
\def\y{1}
\def\x{1}
\begin{array}{cc}
\begin{array}{c}
\begin{tikzpicture}[every node/.style={outer sep=1,inner sep=0}]
\node (t) at (0,{2*\y}){\small{$t$}};
\node (x) at (0,{1*\y}){\small{$\hat x $}};
\node (f) at (0,{0*\y}){\small{$\chi$}};
\path[decoration={markings,mark=at position 0.5 with {\arrow{>}}}]
(t) edge[postaction={decorate}] (x)
(x) edge[postaction={decorate}] (f)
;
\draw[decoration={markings,mark=at position 0.55 with {\arrow{>}}},postaction={decorate}]
(f.-112) .. controls ({0-0.5},{0-0.7}) and ({{0+0.5}},{0-0.7}) .. (f.-68);
\end{tikzpicture}
\end{array}
&
\begin{array}{c}
\begin{tikzpicture}[every node/.style={outer sep=1,inner sep=0}]
\node (ff) at ({-1.5*\x},{2*\y}){\small{$\chi\chi$}};
\node (fb) at ({-0.5*\x},{2*\y}){\small{$\chi\hat x $}};
\node (bf) at ({+0.5*\x},{2*\y}){\small{$\hat x \chi$}};
\node (bb) at ({+1.5*\x},{2*\y}){\small{$\hat x \hat x $}};
\node (ft) at ({+0.0*\x},{1*\y}){\small{$\chi t$}};
\node (bt) at ({+0.0*\x},{0*\y}){\small{$\hat x  t$}};
\path[decoration={markings,mark=at position 0.5 with {\arrow{>}}}]
(ff) edge[postaction={decorate}] (ft)
(fb) edge[postaction={decorate}] (ft)
(bf) edge[postaction={decorate}] (ft)
(bb) edge[postaction={decorate}] (ft)
(ft) edge[postaction={decorate}] (bt)
;
\draw[decoration={markings,mark=at position 0.55 with {\arrow{>}}},postaction={decorate}]
(bt.-112) .. controls ({0-0.5},{0-0.7}) and ({{0+0.5}},{0-0.7}) .. (bt.-68);
\end{tikzpicture}
\end{array}
\\
\X&B
\end{array}
\]
{\caption{\label{fig:B} Illustration of $B$ (for $X$ as in Figure \ref{fig:A}).}}
\end{figure}

%% file: figY.tex
\begin{figure}[p]
\[
\begin{array}{c}
\def\y{1.6}
\def\x{1.2}
\begin{array}{c}
\begin{tikzpicture}[every node/.style={outer sep=1,inner sep=0}]
\node (0ff) at ({-4.5*\x},{2*\y}){\small{$0,\chi\chi$}};
\node (1ff) at ({-3.5*\x},{2*\y}){\small{$1,\chi\chi$}};
\node (0fb) at ({-2.5*\x},{2*\y}){\small{$0,\chi\hat x $}};
\node (1fb) at ({-1.5*\x},{2*\y}){\small{$1,\chi\hat x $}};
\node (0bf) at ({-0.5*\x},{2*\y}){\small{$0,\hat x \chi$}};
\node (1bf) at ({+0.5*\x},{2*\y}){\small{$1,\hat x \chi$}};
\node (0bb) at ({+1.5*\x},{2*\y}){\small{$0,\hat x \hat x $}};
\node (1bb) at ({+2.5*\x},{2*\y}){\small{$1,\hat x \hat x $}};
\node (ff) 	at ({+3.5*\x},{2*\y}){\small{\textcolor{blue}{$\chi \chi$}}};
\node (fb) 	at ({+4.5*\x},{2*\y}){\small{\textcolor{blue}{$\chi \hat x $}}};
\node (1ft) at ({+0.5*\x},{1*\y}){\small{$1,\chi t$}};
\node (0ft) at ({-0.5*\x},{1*\y}){\small{$0,\chi t$}};
\node (1bt) at ({-1.5*\x},{1*\y}){\small{$1,\hat x  t$}};
\node (ft) 	at ({+1.5*\x},{1*\y}){\small{\textcolor{blue}{$\chi t$}}};
\node (0bt) at ({+0.0*\x},{0*\y}){\small{$0,\hat x  t$}};
\path[decoration={markings,mark=at position 0.5 with {\arrow{>}}}]
(0ff) edge[postaction={decorate}] (0ft) 
(1ff) edge[postaction={decorate}] (0ft)
(0fb) edge[postaction={decorate}] (0ft)
(1fb) edge[postaction={decorate}] (0ft)
(0bf) edge[postaction={decorate}] (0ft)
(1bf) edge[postaction={decorate}] (0ft)
(0bb) edge[postaction={decorate}] (0ft)
(1bb) edge[postaction={decorate}] (0ft)
(0ft) edge[postaction={decorate}] (0bt)
(1bt) edge[postaction={decorate}] (0bt)
(1ft) edge[postaction={decorate}] (0bt)
(ff)  edge[postaction={decorate},blue] (0ft)
(fb)  edge[postaction={decorate},blue] (0ft)
(ft)  edge[postaction={decorate},blue] (0bt)
;
\draw[decoration={markings,mark=at position 0.55 with {\arrow{>}}},postaction={decorate}] 
(0bt.-112) .. controls ({0-0.5},{0-0.7}) and ({{0+0.5}},{0-0.7}) .. (0bt.-68);
\end{tikzpicture}
\end{array}
\\
Y
\end{array}
\]
{\caption{\label{fig:Y} Illustration of $Y$ (for $X$ as in Figure \ref{fig:A}), obtained from $PB$ (in black) by adding the blue vertices (that of $V_\chi(B)$) and the blue edges.}}
\end{figure}
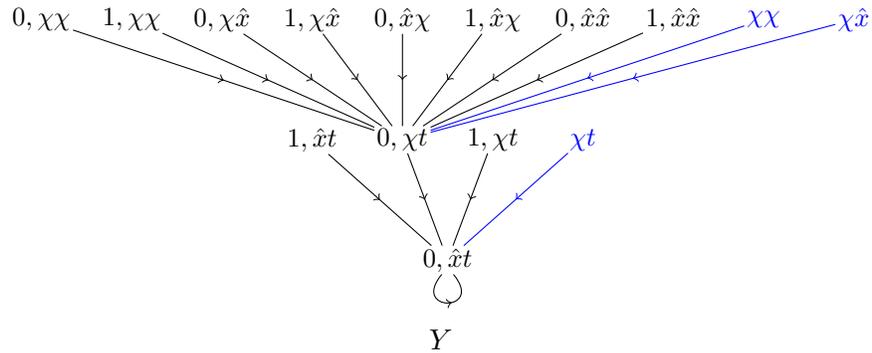

%% file: figXYAB.tex
\begin{figure}[p]
\[
\begin{array}{c}
\def\y{2.1}
\def\x{1.5}
\def\xx{1.75}
\begin{array}{c}
\begin{tikzpicture}[every node/.style={outer sep=0,inner sep=2},every text node part/.style={align=center}]
\node[draw,rectangle,rounded corners] (0ff) at ({-4.5*\x},{2*\y})
{\scriptsize{$\hat x,(0,\chi\chi)$}\\[1mm]\scriptsize{$\chi,(0,\chi\chi)$}};
\node[draw,rectangle,rounded corners] (1ff) at ({-3.5*\x},{2*\y})
{\scriptsize{$\hat x,(1,\chi\chi)$}\\[1mm]\scriptsize{$\chi,(1,\chi\chi)$}};
\node[draw,rectangle,rounded corners] (0fb) at ({-2.5*\x},{2*\y})
{\scriptsize{$\hat x,(0,\chi\hat x)$}\\[1mm]\scriptsize{$\chi,(0,\chi\hat x)$}};
\node[draw,rectangle,rounded corners] (1fb) at ({-1.5*\x},{2*\y})
{\scriptsize{$\hat x,(1,\chi\hat x)$}\\[1mm]\scriptsize{$\chi,(1,\chi\hat x)$}};
\node[draw,rectangle,rounded corners] (0bf) at ({-0.5*\x},{2*\y})
{\scriptsize{$\hat x,(0,\hat x \chi)$}\\[1mm]\scriptsize{$\chi,(0,\hat x \chi)$}};
\node[draw,rectangle,rounded corners] (1bf) at ({+0.5*\x},{2*\y})
{\scriptsize{$\hat x,(1,\hat x \chi)$}\\[1mm]\scriptsize{$\chi,(1,\hat x \chi)$}};
\node[draw,rectangle,rounded corners] (0bb) at ({+1.5*\x},{2*\y})
{\scriptsize{$\hat x,(0,\hat x \hat x)$}\\[1mm]\scriptsize{$\chi,(0,\hat x \hat x)$}};
\node[draw,rectangle,rounded corners] (1bb) at ({+2.5*\x},{2*\y})
{\scriptsize{$\hat x,(1,\hat x \hat x)$}\\[1mm] \scriptsize{$\chi,(1,\hat x \hat x)$}};
\node[draw,rectangle,rounded corners] (ff) 	at ({+3.5*\x},{2*\y})
{\scriptsize{$\hat x,\chi \chi$}\\[1mm]\scriptsize{$\chi,\chi \chi$}};
\node[draw,rectangle,rounded corners] (fb)	at ({+4.5*\x},{2*\y})
{\scriptsize{$\hat x,\chi \hat x$}\\[1mm]\scriptsize{$\chi,\chi \hat x$}};
\node[draw,rectangle,rounded corners] (1ft) at ({+0.5*\x*\xx},{1*\y})
{\scriptsize{$\hat x,(1,\chi t)$}\\[1mm]\scriptsize{$\chi,(1,\chi t)$}};
\node (x0ft) at ({-1.5*\x*\xx},{1*\y})
{\scriptsize{$\hat x,(0,\chi t)$}};
\node (0ft) at ({-0.5*\x*\xx},{1*\y})
{\scriptsize{$\chi,(0,\chi t)$}};
\node[draw,rectangle,rounded corners] (1bt) at ({-2.5*\x*\xx},{1*\y})
{\scriptsize{$\hat x,(1,\hat x  t)$}\\[1mm]\scriptsize{$\chi,(1,\hat x  t)$}};
\node [draw,rectangle,rounded corners] (ft)	at ({+1.5*\x*\xx},{1*\y})
{\scriptsize{$\hat x,\chi t$}\\[1mm]\scriptsize{$\chi,\chi t$}};
\node (x0bt) at ({+2.5*\x*\xx},{1*\y})
{\scriptsize{$\hat x,(0,\hat x  t)$}};
\node (0bt) at ({+0.0*\x},{0*\y}){\scriptsize{$\chi,(0,\hat x  t)$}};
\path[decoration={markings,mark=at position 0.5 with {\arrow{>}}}]
(0ff.-90) 	edge[postaction={decorate}] (0ft) 
(1ff.-90) 	edge[postaction={decorate}] (0ft)
(0fb.-90) 	edge[postaction={decorate}] (0ft)
(1fb.-90) 	edge[postaction={decorate}] (0ft)
(0bf.-90) 	edge[postaction={decorate}] (0ft)
(1bf.-90) 	edge[postaction={decorate}] (0ft)
(0bb.-90) 	edge[postaction={decorate}] (0ft)
(1bb.-90) 	edge[postaction={decorate}] (0ft)
(0ft.-90) 	edge[postaction={decorate}] (0bt)
(1bt) 		edge[postaction={decorate}] (0bt)
(1ft) 		edge[postaction={decorate}] (0bt)
(ff.-90)  	edge[postaction={decorate}] (0ft)
(fb.-90)  	edge[postaction={decorate}] (0ft)
(ft)  		edge[postaction={decorate}] (0bt)
(x0ft) 		edge[postaction={decorate}] (0bt)
(x0bt) 	edge[postaction={decorate}] (0bt)
;
\draw[decoration={markings,mark=at position 0.55 with {\arrow{>}}},postaction={decorate}] 
(0bt.-112) .. controls ({0-0.5},{0-0.7}) and ({{0+0.5}},{0-0.7}) .. (0bt.-68);
\end{tikzpicture}
\end{array}
\\
XY
\end{array}
\]
\[
\begin{array}{c}
\def\y{2.1}
\def\x{1.5}
\def\xx{1.75}
\begin{array}{c}
\begin{tikzpicture}[every node/.style={outer sep=0,inner sep=2},every text node part/.style={align=center}]
\node[draw,rectangle,rounded corners] (0ff) at ({-4.5*\x},{2*\y})
{\scriptsize{$(\hat x,0),\chi\chi$}\\[1mm]\scriptsize{$(\chi,0),\chi\chi$}};
\node[draw,rectangle,rounded corners] (1ff) at ({-3.5*\x},{2*\y})
{\scriptsize{$(\chi,1),\hat x\chi$}\\[1mm]\scriptsize{$(\chi,1),\chi\chi$}};
\node[draw,rectangle,rounded corners] (0fb) at ({-2.5*\x},{2*\y})
{\scriptsize{$(\hat x,0),\chi\hat x$}\\[1mm]\scriptsize{$(\chi,0),\chi\hat x$}};
\node[draw,rectangle,rounded corners] (1fb) at ({-1.5*\x},{2*\y})
{\scriptsize{$(\chi,1),\hat x\hat x$}\\[1mm]\scriptsize{$(\chi,1),\chi\hat x$}};
\node[draw,rectangle,rounded corners] (0bf) at ({-0.5*\x},{2*\y})
{\scriptsize{$(\hat x,0),\hat x \chi$}\\[1mm]\scriptsize{$(\chi,0),\hat x \chi$}};
\node[draw,rectangle,rounded corners] (1bf) at ({+0.5*\x},{2*\y})
{\scriptsize{$(\hat x,1),\hat x \chi$}\\[1mm]\scriptsize{$(\hat x,1),\chi\chi$}};
\node[draw,rectangle,rounded corners] (0bb) at ({+1.5*\x},{2*\y})
{\scriptsize{$(\hat x,0),\hat x \hat x$}\\[1mm]\scriptsize{$(\chi,0),\hat x \hat x$}};
\node[draw,rectangle,rounded corners] (1bb) at ({+2.5*\x},{2*\y})
{\scriptsize{$(\hat x,1),\hat x \hat x$}\\[1mm] \scriptsize{$(\hat x,1),\chi\hat x$}};
\node[draw,rectangle,rounded corners] (ff) 	at ({+3.5*\x},{2*\y})
{\scriptsize{$u,\hat x \chi$}\\[1mm]\scriptsize{$u,\chi \chi$}};
\node[draw,rectangle,rounded corners] (fb)	at ({+4.5*\x},{2*\y})
{\scriptsize{$u,\hat x\hat x$}\\[1mm]\scriptsize{$u,\chi \hat x$}};
\node[draw,rectangle,rounded corners] (1ft) at ({+0.5*\x*\xx},{1*\y})
{\scriptsize{$(\chi,1),\hat x t$}\\[1mm]\scriptsize{$(\chi,1),\chi t$}};
\node (x0ft) at ({-1.5*\x*\xx},{1*\y})
{\scriptsize{$(\hat x,0),\chi t$}};
\node (0ft) at ({-0.5*\x*\xx},{1*\y})
{\scriptsize{$(\chi,0),\chi t$}};
\node[draw,rectangle,rounded corners] (1bt) at ({-2.5*\x*\xx},{1*\y})
{\scriptsize{$(\hat x,1),\hat x  t$}\\[1mm]\scriptsize{$(\hat x,1),\chi t$}};
\node [draw,rectangle,rounded corners] (ft)	at ({+1.5*\x*\xx},{1*\y})
{\scriptsize{$u,\hat x t$}\\[1mm]\scriptsize{$u,\chi t$}};
\node (x0bt) at ({+2.5*\x*\xx},{1*\y})
{\scriptsize{$(\hat x,0),\hat x  t$}};
\node (0bt) at ({+0.0*\x},{0*\y}){\scriptsize{$(\chi,0),\hat x  t$}};
\path[decoration={markings,mark=at position 0.5 with {\arrow{>}}}]
(0ff.-90) 	edge[postaction={decorate}] (0ft) 
(1ff.-90) 	edge[postaction={decorate}] (0ft)
(0fb.-90) 	edge[postaction={decorate}] (0ft)
(1fb.-90) 	edge[postaction={decorate}] (0ft)
(0bf.-90) 	edge[postaction={decorate}] (0ft)
(1bf.-90) 	edge[postaction={decorate}] (0ft)
(0bb.-90) 	edge[postaction={decorate}] (0ft)
(1bb.-90) 	edge[postaction={decorate}] (0ft)
(0ft.-90) 	edge[postaction={decorate}] (0bt)
(1bt) 		edge[postaction={decorate}] (0bt)
(1ft) 		edge[postaction={decorate}] (0bt)
(ff.-90)  	edge[postaction={decorate}] (0ft)
(fb.-90)  	edge[postaction={decorate}] (0ft)
(ft)  		edge[postaction={decorate}] (0bt)
(x0ft) 		edge[postaction={decorate}] (0bt)
(x0bt) 	edge[postaction={decorate}] (0bt)
;
\draw[decoration={markings,mark=at position 0.55 with {\arrow{>}}},postaction={decorate}] 
(0bt.-112) .. controls ({0-0.5},{0-0.7}) and ({{0+0.5}},{0-0.7}) .. (0bt.-68);
\end{tikzpicture}
\end{array}
\\
AB
\end{array}
\]
{\caption{\label{fig:XYAB} Illustration of $XY$ and $AB$ (for $X$ as in Figure \ref{fig:A}). Vertices in the same rectangle have the same out-neighbor. The drawing of $AB$ is obtained from that of $XY$ by relabeling  the vertices of $XY$ with the function $\phi$.}}
\end{figure}
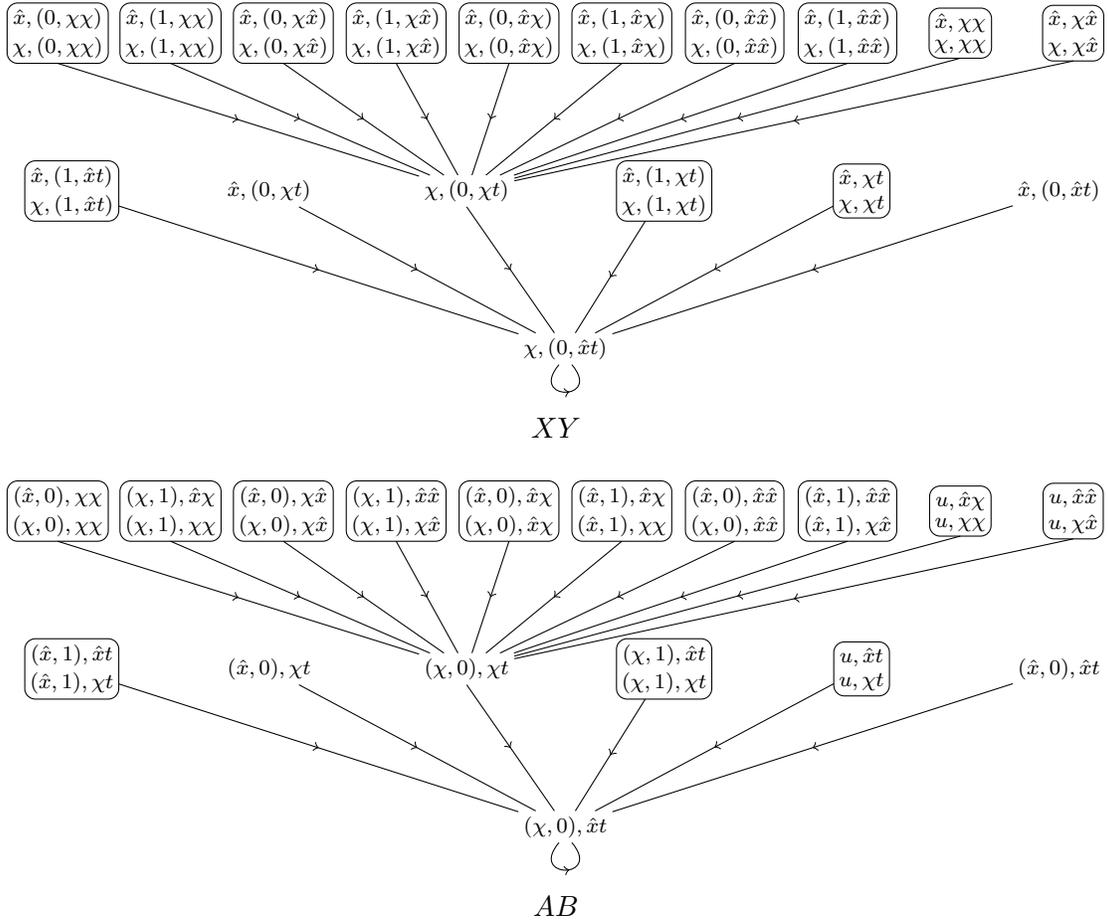